\documentclass{article}
\usepackage{amsmath, epsfig,amssymb, multicol}
\usepackage{color}
\newtheorem{lemma}{Lemma}
\newtheorem{theorem}{Theorem}

\renewcommand{\Pr}{{\rm Pr}}
\newcommand{\E}{{\rm E}}
\renewcommand{\l}{{\cal L}}
\newcommand{\A}{{\cal A}}
\newcommand{\ignore}[1]{}
\renewcommand{\ln}{\log}

\title{Decomposition of random graphs into complete bipartite graphs}
\author{Fan Chung \thanks{Department of Mathematics, University of California, San Diego,  La Jolla, CA 92093,
({\tt fan@ucsd.edu}), Research is supported in part by ONR MURI N000140810747, and AFSOR AF/SUB 552082.}
\and Xing Peng  \thanks{Center for Applied Mathematics, Tianjin University,  Tianjin,    300072 China,
({\tt x2peng@tju.edu.cn}).}
}

\date{}
\begin{document}
\maketitle

\begin{abstract}
We consider the problem of partitioning  the edge set of a graph $G$
into the minimum number   $\tau(G)$ of   edge-disjoint complete bipartite subgraphs.  We show that for a random graph $G$ in $G(n,p)$, where $p$ is a  constant no greater than $1/2$, asymptotically almost surely $\tau(G)$ is between $n- c(\ln_{1/p} n)^{3+\epsilon}$ and
$n - (2+o(1))\ln_{1/(1-p)} n$  for any positive constants $c$ and $\epsilon$.

\end{abstract}

\section{Introduction}
For a graph $G$, the {\it bipartition number}, denoted by $\tau(G)$,  is the minimum number of complete bipartite subgraphs that are edge-disjoint and whose union is    the edge set of $G$. In 1971,  Graham and Pollak \cite{gp1}  proved that
\begin{eqnarray}
\label{gp0}
\tau(K_n) = n-1.
\end{eqnarray}
In particular, they showed that for a graph $G$ on $n$ vertices, the bipartition number $\tau(G)$ is bounded below as follows:
\begin{eqnarray}
\label{gp}
\tau(G) \geq \max \left\{n_+, n_-\right\}
\end{eqnarray}
where $n_+$ is  the number of positive eigenvalues and  $n_-$ is the number of negative eigenvalues of the adjacency matrix of $G$.
Then, (\ref{gp0}) follows from (\ref{gp}).
Since then, there have been a number of alternative proofs for (\ref{gp0}) by using linear algebra \cite{ lo, peck,tv} or by using matrix enumeration \cite{v1,v2}.

Let $\alpha(G)$ denote the independence number of $G$, which is the maximum number of vertices so that there are no edges among some set of  $\alpha(G)$ vertices in $G$.
A star is a special bipartite graph in which all edges share a common vertex which we call the center of the star.  For a graph $G$ on $n$ vertices, the edge set of $G$ can obviously be
decomposed into $n- \alpha(G)$ stars centered at vertices in the complement of a largest independent set. It follows immediately that
\begin{eqnarray}
\label{ep}
\tau(G) \leq n -\alpha(G).
 \end{eqnarray}

It was mentioned in a 1988 paper \cite{krw}, by Kratzke, Reznick and West,  that Erd\H{o}s conjectured  the equality in (\ref{ep}) holds for almost all graphs $G \in G(n,1/2)$, although we could not find this conjecture in any other publication on Erd\H{o}s' problems.

%Erd\H{o}s conjectured   (see \cite{krw}) that   for a random graph $G$ in $G(n,1/2)$,   the equality in (\ref{ep})  asymptotically almost surely holds:
%
%\vspace{.1in}
%
%\noindent
%{\bf A conjecture  of Erd\H{o}s:}~~
%{\it Almost all graphs $G$ on $n$ vertices satisfy}
%\begin{eqnarray} \label{ep1} \tau(G) = n - \alpha(G). \end{eqnarray}
%We remark that the above conjecture by Erd\H{o}s implies that a decomposition by stars is the  best possible for almost all graphs.

Let $\beta(G)$ be the size of the largest induced complete bipartite graph in $G$.  Another possible upper bound for $\tau(G)$, as pointed out by Alon \cite{alon}, is
\[
\tau(G) \leq n-\beta(G)+1.
\]
This follows from the fact that  edges in $G$  can be partitioned into $n-\beta(G)$ stars and a largest induced complete bipartite graph. Therefore we get
\[
\tau(G) \leq \min \{n-\alpha(G),n-\beta(G)+1\}.
\]
A random graph  asymptotically almost surely has  an independent set of order $c \ln n$ and therefore $\tau(G) \leq n - c \ln n$.
For the lower bound,
 for a random graph $G$, it is well known  that asymptotically almost surely  the number of positive and negative eigenvalues
 of the adjacency matrix of $G$ is bounded above by $n/2+ c' \sqrt n$. Consequently,  the inequality in (\ref{gp})  yields  a rather weak  lower bound of $\tau(G)$.
We will prove the following theorem.

\begin{theorem} \label{t:thm0}
For a random graph $G$ in $G(n,1/2)$,   asymptotically almost surely  the bipartition number $\tau(G)$ of $G$ satisfies
\[
n-c(\ln_2 n)^{3+\epsilon} \leq \tau(G)\leq n- (2+o(1)) \ln_{2} n
\]
 for  any positive constants $c$ and $\epsilon$.
\end{theorem}
Theorem \ref{t:thm0} is a special case of the following theorem.

\begin{theorem} \label{t:thm1}
For a random graph $G$ in $G(n,p)$, where   $p$ is a  constant no greater than $1/2$ ,  asymptotically almost surely  the bipartition number $\tau(G)$ of $G$ satisfies
\[
n- c(\ln_{1/p} n)^{3+\epsilon} \leq \tau(G) \leq n-(2+o(1)) \ln_{1/(1-p)} n
\]
 for  any positive constants  $c$ and $\epsilon$.
\end{theorem}

We remark that our techniques can be extended to the case where $p \leq 1-c$ for any positive constant $c$, but we restrict our attention here to the case where $p \leq 1/2$.

 Alon \cite{alon} disproved Erd\H{o}s' conjecture by showing asymptotically almost surely $\tau(G) \leq n-\alpha(G)-1$  for most values of $n$ if $G \in G(n,1/2)$. Recently, Alon, Bohman, and Huang \cite{abh} established a better upper bound which asserts if $G \in G(n,1/2)$, then asymptotically almost surely  $\tau(G) \leq n-(1+c) \alpha(G)$ for some  small positive constant $c$. This result   implies that Erd\H{o}s' conjecture is even false.
For sparser random graphs,
 Alon \cite{alon}    proved   that there exists some (small) constant $c$ such that for $\tfrac{2}{n} \leq p \leq c$, the bipartition number for a random graph $G$ in $G(n,p)$ satisfies $\tau(G)=n-\Theta\left(\tfrac{\log np}{p}  \right)$ asymptotically almost surely.

 We remark that the difficulty for computing $\tau(G)$  is closely related to the intractability of computing $\alpha(G)$. In general, the problem of determining $\alpha(G)$ is an
NP-complete problem, as one of the original 21 NP-complete problems in Karp \cite{karp}. If $G$ does not contain
a $4$-cycle, then $\tau(G)=n-\alpha(G)$. Schrijver showed that  the problem of determining $\alpha(G)$ for the family of $C_4$-free graphs $G$ remains  NP-complete \cite{krw}. Therefore the problem of determining $\tau(G)$ is also NP-complete.
Nevertheless, Theorem \ref{t:thm1}  implies that for almost all graphs $G$, we can bound $\tau(G)$ within a relatively small range.

We also consider a variation of the bipartition number by requiring an additional condition that
no complete  bipartite graph in the partition is  a star. We define the {\it strong bipartition number}, denoted by $\tau'(G)$,  to be the minimum number of complete bipartite graphs (which are not stars) needed to partition the edge set of $G$.  It is possible that a graph $G$ dose not admit  such a partition, then we define $\tau'(G)$ as $\infty$ in this case; if $|V(G)| \leq 2$, then we define $\tau'(G)$ to be zero. We will show that for a random graph $G \in G(n,p)$,  the strong bipartition number satisfies $\tau'(G) \geq 1.0001 n$ if $p$ is a constant and $p \leq \tfrac{1}{2}$.

The paper is organized as follows: In the next section, we state some  definitions and basic facts that we will use later.
In Section \ref{sec3}, we establish  upper bounds for  the number of edges covered by several specified families of complete bipartite subgraphs. In Section \ref{sec4}, we consider the remaining uncovered edges  and give corresponding lower bounds that our main theorem needs. In Section \ref{sec5}, we show that asymptotically almost surely  the strong bipartition number is at least $1.0001 n$ for a random graph on $n$ vertices. In Section \ref{sec6}, we use  the  lemmas and the strong bipartition theorem to prove Theorem \ref{t:thm1}.   A number of problems and remarks are mentioned in Section
\ref{sec8}.
% In the appendix, a combinatorial proof for the Graham-Pollak Theorem in (\ref{gp0}) is given.

\section{Preliminaries}
  Let $G=(V,E)$ be a graph. For a vertex $v \in V(G)$,  the neighborhood $N_G(v)$ of  $v$ is  the set $\{u \colon u \in V(G) \textrm{ and } \{u,v\} \in E(G)\}$ and the degree $d_G(v)$ of $v$ is $|N_G(v)|$.  For a hypergraph $H=(V,E)$ and $v \in V(H)$,  we define the degree $d_H(v)$ to be  $|\{F \colon v \in F \textrm{ and } F \in E(H) \}|$.
  For $U \subseteq V(G)$, let  $e(U)$  be the number of edges of $G$ with both endpoints in $U$ and  $G[U]$ be the subgraph induced by $U$. Furthermore,  $2^{U}$ denotes the power set of $U$.
For two subsets  $A$ and $B$  of $V$, we define $E(A,B)=\{ \{u,v\} \in E \colon u\in A \textrm{ and } v \in B\}$.  We say $A$ and $B$ form a complete bipartite graph if $A \cap B=\emptyset$ and $\{u,v\} \in E(G)$ for all $u \in A$ and $v \in B$.

 We will  use the following versions of  Chernoff's inequality and Azuma's  inequality.
\begin{theorem}{\cite{chernoff}}
\label{t:chernoff}
 Let $X_1,\ldots,X_n$ be independent random variables with
$$\Pr(X_i=1)=p_i, \qquad \Pr(X_i=0)=1-p_i.$$
We consider the sum $X=\sum_{i=1}^n X_i$
with expectation $\E(X)=\sum_{i=1}^n p_i$. Then we have
\begin{eqnarray*}
\mbox{(Lower tail)~~~~~~~~~~~~~~~~~}
\qquad \qquad  \Pr(X \leq \E(X)-\lambda)&\leq& e^{-\lambda^2/2\E(X)},\\
\mbox{(Upper tail)~~~~~~~~~~~~~~~~~}
\qquad \qquad
\Pr(X \geq \E(X)+\lambda)&\leq& e^{-\frac{\lambda^2}{2(\E(X) + \lambda/3)}}.
\end{eqnarray*}
\end{theorem}

\begin{theorem} \cite{azuma}
 \label{t:azuma}
Let $X$ be a random variable determined by $m$ trials $T_1,\ldots,T_m$, such that for each $i$, and any two possible sequences of outcomes $t_1,\ldots,t_i$ and $t_1,\ldots,t_{i-1}, t_i'$:
\[
|\E\left(X|T_1=t_1,\ldots,T_i=t_i\right) -\E\left(X|T_1=t_1,\ldots, T_{i-1}=t_{i-1},T_i=t_i'\right)   | \leq c_i
\]
then
\[
\Pr\left(|X-\E(X)| \geq \lambda  \right) \leq 2 e^{-{\lambda}^2/2\sum_{i=1}^m c_i^2}.
\]
\end{theorem}

The following lemma on edge density will be useful later.
\begin{lemma} \label{l:lm1}
Asymptotically almost surely a   random graph $G$ in  $G(n,p)$ satisfies, for all  $U \subset V(G)$ with $|U| \geq \sqrt{ \ln n}$,
\[
\left|e(U)-\frac{p}{2}|U|^2 \right| \leq C |U|^{3/2} \ln^{1/2} n
\]
 where $C$ is some positive constant.
% with $|U| \geq n^{0.67}$.
\end{lemma}
The lemma follows from  Theorem \ref{t:chernoff}.

The following lemma is along the lines of a classical result of Erd\H{o}s for random graphs \cite{ep}. We include the statement and a short proof here for the sake of completeness.

\begin{lemma} \label{l:lm2}
For  $G \in G(n,p)$, where  $p$ is a constant no greater than $1/2$,
asymptotically almost surely all complete bipartite graphs $K_{A,B}$ in $G$ with $|A| \leq |B|$  satisfy  $|A| \leq 2 \ln_b n$ where $b=1/p$.
\end{lemma}

\noindent
{\bf Proof:}
%It suffices to prove for the case $p=\tfrac{1}{2}$.
For  two subsets $A$ and $B$ of $V(G)$,  with $|A|=|B| =k$,
the probability that
$A$ and $B$ form a complete bipartite graph in $G(n,p)$  is at most
$p^{k^2}.$
There are at most
$ \binom{n}{k} \binom{n}{ k}$ choices for $A$ and $B$.
For $k \geq 2 \ln_b n$, we have
\[
\binom{n}{k}^2p^{k^2}=o(1) \]
as $p$ is a constant.
The lemma then follows.
 \hfill $\square$

%We note that Erd\H{o}s'  result on
%Ramsey's theorem \cite{ep} states that every 2-coloring of the edges of the complete graph $K_n$ contains a monochromatic clique of order $\tfrac{1}{2}\ln_2 n$.  It is not hard to show along the same lines that  $G(n,1/2)$ contains a complete bipartite graph $K_{A,B}$ with $|A|=|B|=\tfrac{1}{4} \ln_2 n$   and the bound in Lemma \ref{l:lm2} is tight up to a constant factor.

The upper bound in Theorem \ref{t:thm1} is an immediate consequence  of (\ref{ep}).
%and the classical results on the independence number  $\alpha(G)$ of a random graph.
%An upper bound for $\alpha(G)$ can be found in \cite{ep}.
   The problem of determining the independence number for a random graph has been extensively studied in the literature.
The asymptotic order of $\alpha(G)$ for $G $ in $G(n,p)$  was determined in \cite{gm}.

\begin{theorem}\label{alpha} \cite{gm}
 If  $p$ is a constant, $p < 1-c$, and $G \in G(n,p)$, then  asymptotically almost surely $\alpha(G)$ is of order
\[ \alpha(G)= 2 \log_{1/(1-p)} n + o(\log n) \]
where  $c$ is a positive constant .
%\item[(2) \cite{frieze}] If $p=o(1)$, $\alpha(G)$ asymptotically almost surely satisfies
%  \[ \alpha(G)= \frac 2 p \log np -  (1+o(1)) \log \log np .\]
%  where $\log$ denotes the natural logarithm.
%\end{description}
\end{theorem}

\section{Edges covered by a given family of subsets}\label{sec3}
For a graph $G=(V,E)$  and $A \subset V$, we define
 \[
 V(G,A)=\{v:  v \in V(G) \setminus A \textrm{ and } \{u, v\} \in E  \textrm{ for all } u \in A \}.
 \]
 It immediately follows    that  $A$ and $B$ form a complete bipartite graph if $B$ is contained in $V(G,A)$, namely,  $B \subseteq V(G,A)$.
 We say an edge $\{u,v\} \in E$ is covered by $A$ if either  $u \in A$ and $v \in V(G,A)$ or $v \in A$ and $u \in V(G,A)$.

\begin{figure}[ht\tau]
 \centerline{\psfig{figure=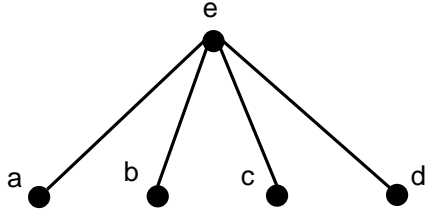, width=0.4\textwidth} }
\caption{An illustration of $f(G,\A)$.}
\label{fg:fg1}
\end{figure}

 For  $\mathcal A=\{A_1,A_2,\ldots,A_k\} \subseteq 2^V$ and  $\sigma$,  a linear ordering of $[k]$,  we define a function $l$ as  follows.  For notational convenience, we  use $i$ to denote the $i$-th element under the ordering $\sigma$.   For each $1 \leq i \leq k$,
 we define  $G_{i}$ and $l(i)$ recursively.
We  let $G_{1}=G$ and let $l(1)$ be   an arbitrary  subset of $V(G_{1},A_{1})$. Given $G_{i-1}$, we let $G_{i}$ be a new graph with the vertex set $V(G)$ and the edge set $E(G_{i-1}) \setminus E(A_{i-1}, l(i-1))$. We set $l(i)$ to be an arbitrary subset of $V(G_{i},A_{i})$. We define
\[
f(G,\A)= \underset{\sigma}   {\max} \ \underset{l} {\max} \ \sum_{i=1}^k |E(A_i,l(i))|.
\]
Basically, for given $A_i$'s, we wish to choose $B_i$'s so that the   complete bipartite graphs $K_{A_1,B_1}, \ldots,K_{A_k,B_k}$ cover as many  edges in $G$ as possible.
An example is illustrated in   Figure \ref{fg:fg1} for  $\A=\{A_1,A_2,A_3\}$ with
$A_1=\{a,b\},A_2=\{b,c\},$ and $A_3=\{c,d\}$.  Here $f(G,\A)=4$ is achieved by $\sigma=\textrm{identity}$, $l(1)=\{e\}$, $l(2)=\emptyset$, and $l(3)=\{e\}$, or $\sigma=(213)$, $l(1)=\emptyset$ and  $l(2)=l(3)=\{e\}$.
We observe $f(G,\A) \leq \sum_{i=1}^k |E(A_i,V(G,A_i))|$ as  $l(i) \subseteq V(G,A_i)$.
When $U \subset V(G)$ and $A \subset U$ for each $A \in \A$,  we use $f(G,U,\mathcal A)$ to denote $f(G[U],\A)$. Note that
 $G[U]$ denotes the induced subgraph of $G$ on a subset  $U$ of $V(G)$.

\begin{lemma} \label{l:lm3}
Suppose    $G \in G(n,p)$,   $U \subseteq V(G)$,  and  $\A$ is a family of $2$-sets of $U$ with $|\A| \leq 1.0001|U|$.
If $p \leq \tfrac{1}{2}$, then asymptotically almost surely we have
\[
f(G,U,\A) \leq  2p^2 |\A||U|+  8|U|^{3/2}\log n
\]
for   all choices of $U$ and $\A$.
\end{lemma}

\noindent
{\bf Proof:}  We list  edges with both endpoints in $U$ as $e_1,e_2,\ldots,e_m$ where $m=\binom{|U|}{2}$.
For each  $e_i=\{u_i,v_i\}$ for $1 \leq i \leq m$, we consider  $T_i \in \{\textrm{H},\textrm{T}\}$ where $T_i=\textrm{H}$ means $e_i$ is an edge and $T_i=\textrm{T}$ means $e_i$ is not an edge. To simplify the notation we use $X$ to denote the  random variable $f(G,U,\A)$ and notice that $X$ is determined by $T_1,\ldots,T_m$.  Given the outcome $t_j$ of $T_j$ for each $1 \leq j \leq i-1$ we wish to establish an upper bound for
\begin{equation} \label{eq:lm3}
\left|\E(X|T_1=t_1,\ldots,T_{i-1}=t_{i-1}, T_i=\textrm{H})-\E(X|T_1=t_1,\ldots,T_{i-1}=t_{i-1},T_i=\textrm{T})\right|.
\end{equation}
Let ${\mathcal K}_1$ be the set of graphs over $U$ such that $e_j$ is given by  $t_j$ for each $1 \leq j \leq i-1$ and $e_i$ is a non-edge.
Similarly, let  ${\mathcal K}_2$ be the set of graphs over $U$ such that $e_j$ is given by  $t_j$ for each $1 \leq j \leq i-1$ and $e_i$ is an edge.  We have $|\mathcal K_1|=|\mathcal K_2|$. Thus we get
\[
\E(X|T_1=t_1,\ldots,T_{i-1}=t_{i-1}, T_i=\textrm{T})=\sum_{K \in \mathcal K_1} f(K,\A) \Pr(K \in \mathcal K_1)
\]
and
\[
\E(X|T_1=t_1,\ldots,T_{i-1}=t_{i-1}, T_i=\textrm{H})=\sum_{K \in \mathcal K_2} f(K,\A) \Pr(K \in \mathcal K_2).
\]
Define a mapping $\mu: \mathcal K_1 \to \mathcal K_2$ such that $E(K)$ and $E(\mu(K))$ differ only by $e_i$ for each $K \in \mathcal K_1$. We get $\mu$ is a bijection and $\Pr(K \in \mathcal K_1)=\Pr(\mu(K) \in \mathcal K_2)$. Therefore the expression \eqref{eq:lm3} can be bounded from above by
\[
\sum_{K \in \mathcal K_1} |f(K,\A)-f(\mu(K),\A)| \Pr(K \in \mathcal K_1).
\]
 Notice that each edge can be covered by at most once.  We observe $|f(K,\A)-f(\mu(K),\A)| \leq 2$ because $e_i$ and the other edge sharing one endpoint with $e_i$  could be covered by $\A$ in $\mu(K)$ but not in $K$. Therefore \eqref{eq:lm3} is bounded above by 2.

Now we apply Theorem \ref{t:azuma}  for $\lambda=8|U|^{3/2} \ln n $ and $c_i=2$. Then we have
\begin{equation} \label{eq:azuma}
\Pr\left( \left|X-\E(X) \right| \geq  8|U|^{3/2} \ln n \right) \leq  2e^{-64|U|^{3}\ln^2 n/2\sum_{i=1}^m c_i^2} \leq 2e^{-4|U|\ln^2 n},
\end{equation}
using the fact  $m \leq \tfrac{|U|^2}{2}$. To estimate  $\E(X)$, we note that  $\E(f(G,U,A)) \leq 2p^2 |U| $ for a fixed $A \in \A$.
 %and for each $w \in U \setminus A$ let $X_w$ be the event that $(u,w), (v,w) \in E(G(n,p))$ and $I_w$ be the random indicator variable for %$X_w$ . We have $\Pr(X_w)=p^2$ and $f(G(n,p),U,A)=\sum_{w \in U \setminus A} 2I_w $.
Therefore,
\[
 \E(X) \leq \sum_{A \in \A}  \E\left(f(G,U,A)\right) \leq 2p^2|\A||U|.
\]
Thus \eqref{eq:azuma}  implies
\begin{align*}
\Pr\left(X \geq 2p^2|\A||U|+  8|U|^{3/2} \ln n \right) &\leq  \Pr\left(|X- \E(X)| \geq  8|U|^{3/2} \ln n \right)\\
                                                                              &\leq 2 e^{-4|U|\ln^2 n}.
\end{align*}
 Recall the assumptions  $|\A| \leq 1.0001|U|$ and $|A|=2$ for each $A \in \A$. For  fixed sizes of $U$ and $\A$,  the number of choices for $U$ and $\A$  is at most $n^{|U|} |U|^{2|\A|}$ which is less than $n^{3.5|U|} $.   Therefore the probability that there are some $U$ and $\A$ which violate the assertion in the lemma is at most
$ 1.0001 |U| \times n \times n^{3.5|U|} \times 2e^{-4|U| \ln^2 n }   < 2e^{-\ln^2 n} $
for sufficiently large $n$, which completes the proof of this lemma. \hfill $\square$

The following lemmas for other families of sets $\A$ have   proofs which are quite  similar to the proof of Lemma \ref{l:lm3}. We will sketch proofs  here.
\begin{lemma} \label{l:lm31}
 Suppose  $G \in G(n,p)$, $\A$ is a family of subsets of $U \subseteq V(G)$ satisfying $|\A | \leq 1.0001|U|$ and $2 \leq |A| \leq 2 \ln_2 n$ for each $A \in \A$.   If $p \leq \tfrac{1}{2}$, then asymptotically almost surely  we have
\[
 f(G,U,\A) \leq  2p^2 |\A||U|+  8|U|^{3/2} \ln^{3/2} n  \ln^{1/2} |U|
\]
for all choices of $U$ and $\A$ satisfying $|U| \geq \log^2 n$.
\end{lemma}
\noindent
{\bf Proof:}  We will use the Azuma's inequality.  For  $e_i=\{u_i,v_i\}$,  we define $\mathcal K_1$, $\mathcal K_2$, and a bijection $\mu$ similarly. The only difference is that $|f(K,\A)-f(\mu(K),\A) | \leq 2 \log_2 n$ for each $K \in \mathcal K_1$. This is because $e_i$ and at most other $2\log_2 n-1$ edges sharing the same endpoint with $e_i$ could by covered by $\A$ in $\mu(K)$ but not in $K$.

 %For  $e_i=\{u_i,v_i\}$ we consider two sets.
%\[
%\A_1=\{A \in \A: u_i \in A, v_i \not \in A, \textrm{ and } (v_i,v) \textrm{ is not revealed  for at least one } v \in A  \}
%\]
%and
%\[
%\A_2=\{A \in \A: v_i \in A, u_i \not \in A, \textrm{ and } (u_i,v) \textrm{ is not revealed  for at least one } v \in A  \}.
%\]
%We also form
%\[
%B_1=\{v \in A \colon A \in A_1 \textrm{ and } (v_i,v) \textrm{ is not revealed}\}
%\]
%and
%\[
%B_2=\{v \in A \colon A \in A_2 \textrm{ and } (u_i,v) \textrm{ is not revealed}\}.
%\]
%We observe that in each configuration of edges between $v_i$ and $B_1$  the number of edges covered by $\A$ can increase  by at most $2\ln_2 n$ and the same observation holds for $B_2$
Therefore  the corresponding expression for \eqref{eq:lm3} can be upper bounded by $2\log_2 n$.
We can then estimate $X=f(G,U,\A)$ by applying  Theorem \ref{t:azuma}  with
\[
\lambda=8|U|^{3/2} \ln^{3/2} n  \ln^{1/2} |U| \textrm{  and  } c_i=2\ln_2 n.
\]
 This leads to
\begin{eqnarray}\label{p1}
\Pr\left( \left|X-\E(X) \right| \geq  8|U|^{3/2} \ln^{3/2} n  \ln^{1/2} |U| \right)  \leq 2e^{-4|U|\ln  n \log |U| }.
\end{eqnarray}
Since  $\E(X) \leq \sum_{A \in \A} \E(f(G,U,A)) \leq 2p^2|\A||U|$  and   the number of choices for $U$ and $\A$ can be bounded  from above by
\[
1.0001 |U|  \times n^{1+|U|} \times |U|^{2.0003|U|\ln_2 n} \leq e^{2.5|U| \log n \log |U|}.
\]
Here we used the following simple fact $\sum_{1 \leq t \leq k} \binom{s}{t}=(1+o(1))) \binom{s}{k}$ if $k$ is much smaller than $s$.
Therefore, the probability that there are some $U$ and $\A$ which violate the lemma is at most $e^{- |U| \log n \log |U| } \leq e^{-\log^2 n}$, which completes the proof of   the lemma. \hfill $\square$
\begin{lemma} \label{l:lm32}
Suppose $G \in G(n,p)$, $\A$ is a family of subsets of $U \subseteq V(G)$ satisfying $|\A | \leq 1.0001 |U|$ and $3 \leq |A| \leq 2 \ln_2 n$ for each $A \in \A$.  If  $p \leq \tfrac{1}{2}$, then asymptotically almost surely  we have
\[
 f(G,U,\A)  \leq 3p^3 |\A| |U|+ 8|U|^{3/2} \ln^{3/2} n  \ln^{1/2} |U|
 \]
 for all choices of $U$ and $\A$ satisfying $|U| \geq \log^2 n$.
\end{lemma}
\noindent
{\bf Proof:} The proof is similar to that of  Lemma \ref{l:lm31}. The only difference is that we assume  $|A| \geq 3$ and therefore
\[
\E(f(G,U,\A)) \leq \sum_{A \in \A} \E(f(G,U,A)) \leq 3p^3 |\A||U|.
\]
We use  Theorem \ref{t:azuma} in a similar way as in the proof of Lemma \ref{l:lm31} to  complete the proof of  Lemma \ref{l:lm32}. \hfill $\square$
\begin{lemma} \label{l:lm33}
 Suppose  $G \in G(n,p)$, $\A$ is a family of subsets of $U \subseteq V(G)$ satisfying $|\A | \leq |U|^{1+\delta}$ and $ \delta \ln_b |U| \leq |A| \leq  2 \ln_2 n$ for each $A \in \A$, where $b=\frac{1}{p}$ and $\delta$ is a fixed small  positive constant.  If $p \leq \tfrac{1}{2}$, then asymptotically almost surely  we have
\[
 f(G,U, \A) \leq   \delta |\A||U|^{1-\delta} \ln_b |U|+8 |U|^{(3+\delta)/2} \ln^{3/2} n \log^{1/2} |U|
\]
for all choices of $U$ and $\A$ satisfying $|U| \geq \log^2 n$.
\end{lemma}

\noindent
{\bf Proof:} We use the assumptions on $|A|$ to derive
\[
\E(f(G,U,\A)) \leq \sum_{A \in \A} \E(f(G,U,A)) \leq  \delta  \ln_b |U| p^{\delta \ln_b|U|} |\A||U| \leq \delta |\A||U|^{1-\delta} \ln_b |U|.
\]
Then we bound the number of choices for $U$ and $\A$ from above by
\[
|U|^{1+\delta} \times n^{1+|U|} \times  |U|^{2.01|U|^{1+\delta}\ln_2 n}  \leq e^{3.5 |U|^{1+\delta} \log n \log |U|}.
\]
Applying   Theorem \ref{t:azuma} for $\lambda=8 |U|^{(3+\delta)/2} \ln_2 n \sqrt{ \ln |U| \ln_2 n}$ and $c_i=2 \ln_2 n$,  the lemma then follows.  \hfill $\square$

%We remark that we only use Lemma \ref{l:lm33}  when $p$ being a constant for proving the main theorem.
%XXX We have tow main theorems and we only use Lemma 6 for proving Theorem 1.XXX

\section{Bounding uncovered edges}\label{sec4}

In order to prove the bipartite decomposition theorem, we also need to establish lower bounds for the number of  uncovered edges
for a given family $\A$ of subsets.

%To do so, we will define $g$ and $h$ as follows:

First, we will derive  a lower bound on  the number of uncovered edges for  a collection $\A$  of $2$-sets of $V(G)$ provided $G \in G(n,p)$.
Let $S_0$ be the set of $u \in V(G)$ such that $u$ is in only one $A \in \A$.  For $u$ in $S_0$, we denote the only 2-set containing $u$ by $A_u$.
 Our goal is to give a lower bound on the number of uncovered edges with both endpoints in  $S_0$.  To simplify the estimate,  we impose some technical restrictions  and work on a subset $S$ of $S_0$. To do so, we will lose at most a  factor  of $2$  in the lower bound estimate (which is tolerable).
 To form $S$, for each  $A_u=\{u,v\}$  with $u,v \in S_0$,  we delete one of $u$ and $v$ arbitrarily from $S_0$.
 Let $T=\cup_{u \in S} \{A_u \setminus  \{u\}\}$. Clearly  $S$ and $T$ are disjoint.  We  note that $S$ and $T$ are determined by $\A$. Furthermore,  $|S| \geq |T|$.

 Suppose $G \in G(n,p)$.  For $u,v \in S$, let $X_{u,v}$ be the indicator random variable such that $\{u,v\} \in E(G) $,  $\{u,A_v \setminus \{v\}\}\not \in E(G)$,   and    $\{v,A_u \setminus \{u\}\}\not \in E(G)$.
 Then we define
 \[
 g(G,\A,S,T)=\sum_{\{u,v\} \in \binom{S}{2} } X_{u,v}
 \]
We observe that $g(G,\A,S,T)$ indeed gives a lower bound on the number of edges which are not covered by $\A$.  Since  if $X_{u,v}=1$  and $\{u,v\}$ is covered by some $A \in \A$, then either $A=A_u$ or $A=A_v$. The former case can not happen because we assume $\{v, A_u \setminus \{u\}\}$ is not an edge. We have the similar argument for the latter case.

\begin{lemma} \label{l:lm4}
Suppose  $G \in G(n,p)$,  $U \subseteq V(G)$, and $\A \subseteq \binom{U}{2}$ with $|\A| \leq 1.0001 |U|$.   Let $S$ and $T$ be defined as above. If   $p \leq \tfrac{1}{2}$, then asymptotically almost surely  we have
\[
 g(G,\A, S,T) \geq  p^3 \binom{|S|}{2}  - 4 |U|^{3/2} \log n
\]
for all choices of $U$ and $\A$ with $|U| \geq \log^2 n$.
\end{lemma}

\noindent
{\bf Proof:} We sketch the proof here which is  similar to that of Lemma \ref{l:lm3}.
We list edges with endpoints in $S \cup T$ as $e_1,\ldots,e_m$, where $m=\binom{|S \cup T|}{2}$ and $m \leq \binom{|U|}{2}$.  For each  $e_i=\{u_i,v_i\}$ for $1 \leq i \leq m$, we consider  $T_i \in \{\textrm{H},\textrm{T}\}$ where $T_i=\textrm{H}$ means $e_i$ is an edge and $T_i=\textrm{T}$ means $e_i$ is not an edge.  Let  $X$  denote the   random variable $g(G,\A,S,T)$ for $G \in G(n,p)$. We note   that $X$ is determined  by $T_1,\ldots,T_m$. For the fixed  outcome $t_j$ of  $T_j$ for $1 \leq j \leq i-1$, we consider
\begin{equation} \label{eq:g}
\left|\E(X|T_1=t_1,\ldots,T_{i-1}=t_{i-1},T_i=\textrm{ H}) -  \E(X|T_1=t_1,\ldots,T_{i-1}=t_{i-1},T_i=\textrm{ T})   \right|.
\end{equation}
For  $e_i=\{u_i,v_i\}$, if $u_i,v_i \in T$ then the outcome of $T_i$ does not contribute to \eqref{eq:g}.
If  $u_i,v_i \in S$  then the outcome of $T_i$ can change \eqref{eq:g} by at most one depending on whether $e_i$ is covered or not.
 If  $u_i\in S$ and $v_i \in T$  then the outcome of $T_i$ could effect \eqref{eq:g} by at most one. This is because $e_i$ could make another edge $\{u_i,w\}$  covered by the $2$-set $\{v_i,w\}$.   Thus \eqref{eq:g} is bounded above by  one.

 Applying Azuma's  theorem as stated in Theorem  \ref{t:azuma} with $\lambda=4 |U|^{3/2} \log n $ and $c_i=1$, we have
\begin{equation} \label{eq:eq2}
\Pr\left(|X-\E(X)| \geq  4 |U|^{3/2} \log n \right)   \leq 2 e^{-4|U|\ln^2  n},
\end{equation}
using  $m \leq  \binom{|U|}{2}$.  To estimate  $\E(X)$,
%we define
% $X_{u,v}$, for  $u,v \in S$, to be the event $(u,v) \in E(G)$, $(u,A_v \setminus \{v\}) \not \in E(G)$, and    $(v,A_u \setminus \{u\}) \not \in E(G)$. Here $G \in G(n,p)$ and $A_v$(resp.$A_u$) is the only $2$-set containing $v$(resp.$u$).
%
% Let $I_{u,v}$ denote the random indicator variable for $X_{u,v}$.
we note $\Pr(X_{u,v}=1)=p(1-p)^2 \geq p^3$ as $p \leq \tfrac{1}{2}$.  Thus,
\[
\E(X) = \sum_{\{u,v\} \in \binom{S}{2}} \Pr(X_{u,v}=1) \geq  p^3 \binom{|S|}{2},
\]
which implies
\[
\Pr\left(X \leq p^3 \binom{|S|}{2} - 4 |U|^{3/2} \log n \right) \leq \Pr\left(X \leq \E(X)-4 |U|^{3/2} \log n\right) \leq 2e^{-4|U|\ln^2 n}.
\]
We recall $S$ and $T$ are determined by $\A$. For  fixed sizes of $U$ and $\A$, the number of choices for $U$ and $\A$ is at most
\[
n^{|U|} |U|^{2|\A|}.
\]
As $|\A| \leq 1.0001 |U|$,   the probability that  there are some $U$ and $\A$ which violate this lemma is at most
\[
1.0001 |U| \times n^{1+|U|} \times |U|^{2.0003 |U|} \times 2e^{-4|U| \log^2 n} \leq e^{-0.5|U| \log^2 n} < e^{-\log^2 n}.
\]
The lemma is proved. \hfill $\square$

% The above lemma is mainly  for the case that $\A $ consists of $2$-sets. For general $\A$, we will define $h$ as follows.

Next, we wish to establish
 a lower bound on the number of uncovered edges for  general cases of $\A$.

For  $W \subset U \subset V(G)$,
we consider  $L \colon W \to 2^{U \setminus W}$  with the property  $L(w) \cap L(w')= \emptyset$ for $w,w' \in W$.
We define $h(G, U,W,L)$ to be  the number of edges $\{w,w'\}$ in $G$ such that $w,w' \in W$, $\{w,z\} \not \in E(G)$ for each $z \in L(w')$, and $\{w',z'\} \not \in E(G)$ for each $z' \in L(w)$.  We will use the following Lemma  (later we will show  that  $h(G,U,W,L)$  gives a lower bound for the number of  uncovered edges).

\begin{lemma} \label{l:lm5}
   Suppose $G \in G(n,p)$,  $p$ is a constant no greater than $1/2$,  $b=\tfrac{1}{p}$, and $|U| \geq \log^2 n$. Assume   $W \subset U$ satisfies  $|W|=|U|/\ln_b^2|U|$ and  $L$ as defined above satisfies $1 \leq |L(w)| \leq c\ln_b |U|$ for some positive constant $c$. Then
asymptotically almost surely  we have
\[
h(G,U,W,L) \geq c'|U|^{2-2c}/\ln_b^4 |U|- 2 |U|^{3/2} \sqrt{\ln n}
\]
for all choices of $U$, $W$ and $L$, where   $c'$ is some positive constant.
\end{lemma}

\noindent
{\bf Proof:}  For  $u,v \in W$,  let $X_{u,v}$ denote the event that $\{u,v\} \in E(G)$, $\{u,w\} \not \in E(G)$ for each $w \in L(v)$, and $\{v,z\} \not \in E(G)$ for each $z \in L(u)$. Here $G \in G(n,p)$. The indicator random variable for $X_{u,v}$ is written as   $I_{u,v}$. From the definition of $h$, we have
$h(G,U,W,L)=\sum_{u,v \in W} I_{u,v}=Y$ for $G \in G(n,p)$.   Since  $\Pr(I_{u,v}=1) \geq b^{-1-2c \ln_b |U|}$, we have
\[
\E(Y) \geq b^{-1-2c\ln_b |U|}\binom{|W|}{2} \geq c'|U|^{2-2c}/\ln_b^4 |U|,
\]
for some constant $c'$.
From the definition of   $L$, we have $X_{u,v}$ are independent of   one another. By applying the  Chernoff's  bound for the lower tail in Theorem \ref{t:chernoff} with $\lambda=2 |U|^{3/2} \sqrt{\ln n}$, we have
\begin{align*}
\Pr( Y \leq  c'|U|^{2-2c}/\ln_b^4 |U|- 2 |U|^{3/2} \sqrt{\ln n}) & \leq \Pr(Y \leq \E(Y)-2 |U|^{3/2}\sqrt{\ln n}) \\
                                                                &\leq e^{-\tfrac{4|U|^{3} \ln n}{2 \E(Y)}} \leq e^{-2|U| \ln n},
\end{align*}
using the fact that  $\E(Y) \leq |U|^2$. For a given  size  $|U|$ of $U$, it is straightforward to bound the number of choices for $U$, $W$ and $L$  from above by
\[
n^{|U|}|U|^{|U|/\ln_b^2|U|} |U|^{1.1c|U|/\ln_b |U|} < n^{|U|} |U|^{2c|U|/ \ln_b |U|} \leq n^{1.5|U|},
\]
when $n$ is sufficiently large. The probability that there is some $U$, $W$ and $L$ which violate the lemma is at most $n  e^{-2 |U|\ln n} n^{1.5 |U|} < e^{-\log^2 n} $ if $n$ is large enough. This completes the proof of the lemma. \hfill $\square$

\section{A theorem on strong bipartition decompositions}\label{sec5}
 Recall the strong bipartition number $\tau'(G)$ is the minimum number of complete  bipartite graphs whose edges partition the edge set of $G$ and none of them is a star. We fix the constant $b=\tfrac{1}{p}$ and  we will prove the following  theorem.

\begin{theorem} \label{t:thm2}
Suppose   $G \in G(n,p)$,  $U \subseteq V(G)$ is  a vertex subset  with  $|U| \geq  c (\ln_b n)^{3+\epsilon}$
where $ b=1/p$,  $c$ and $\epsilon$ are   positive constants.
For  $p$ being a positive constant  no greater than $1/2$,  asymptotically almost surely  for all $U$ we have
\[
\tau'(G[U]) \geq 1.0001 |U|.
\]

\end{theorem}
The proof of Theorem \ref{t:thm2} is based on several lemmas which we will first prove.
In this section, we may assume that $G \in G(n,p) $ satisfies the statements in all lemmas in the preceding sections.
By Lemma \ref{l:lm1}, the number of edges in $G[U]$ satisfies  $e(G[U])=(\tfrac{p}{2}+o(1))u^2$, here $|U|=u$.
We will prove Theorem \ref{t:thm2}  by contradiction.
%We suppose that there are $m$ complete bipartite graphs $K_{A_1,B_1},\ldots,K_{A_m,B_m}$  such that the edge set of $G$ is a disjoint union of $E(K_{A_i,B_i})$ and
Suppose
\[
E(G[U])=\bigsqcup_{i=1}^m E(K_{A_i,B_i})
\]
and   $m < \left(1+\tfrac{1}{10000}\right) u$ where  `$\sqcup$' denotes the disjoint union.
We assume further $|A_i| \leq |B_i|$ for each $1 \leq i \leq m$.
 Lemma \ref{l:lm2} implies  $|A_i| \leq 2 \ln_b n \leq   2 \ln_2 n$  for each $1 \leq i \leq m$.

 We define $\mathcal L=\{A_i \colon 1 \leq i \leq m\}$.
%  and note that $\mathcal L$ could be a multi-set.
   We consider three subsets of $\mathcal L$  defined as follows:
\begin{align*}
\l_1&=\{A_i \in \l \colon  |A_i|< \delta_1 \ln_b u\}\\
 \l_2&=\{A_i \in \l \colon |A_i| < \delta_2 \ln_b u\} \\
\l_3&=\{A_i \in \l \colon |A_i|=2\},
\end{align*}
where
\[
 \delta_1=\min\left\{\frac{\epsilon}{4(3+\epsilon)},\frac{1}{200}\right\} \textrm{ and } \delta_2=\frac{\delta_1}{10^4}.
\]

We observe that a typical complete bipartite graph $K_{2,B}$ in a random graph $G(n,1/2)$ roughly  contains  the same  the number of edges as  a typical star does. This is the reason why we define the  set $\l_3$.

We have the following lemma.
\begin{lemma} \label{l:lm6}
If $|\l_2| \leq \left( \tfrac{1}{2}+\tfrac{1}{1500} \right)u$, then  we have $|\l_3| \geq \left( \tfrac{1}{2}+\tfrac{1}{2000} \right)u$.
\end{lemma}

\noindent
{\bf Proof:}  We will first prove the following claim:

\noindent
{\bf Claim 1:}  $|\l_3| \geq \left( \tfrac{1}{2}-\tfrac{1}{250} \right)u$.
\\
{\it Proof of Claim 1:}
Suppose  the contrary.   By   Lemma \ref{l:lm3},  the number of edges covered by  $\l_3$  is at most
\[
2p^2 |\l_3|u+8u^{3/2} \log^{1/2} n.
\]
 By   Lemma \ref{l:lm32}, the number of edges covered by  $ \l_2 \setminus \l_3$ (i.e., $|A_i| \geq 3$) is at most
 \[
 3p^3|\l_2 \setminus \l_3|u+8u^{3/2} \log^{3/2} n \log^{1/2} u.
 \]
  Therefore the total number of edges covered by  $ \l_2$ is at most
\[
2p^2|\l_3|u+3p^3|\l_2 \setminus \l_3|u  + 8u^{3/2} \log^{1/2} n   + 8u^{3/2} \log^{3/2} n \log^{1/2} u.
\]
Thus  the number of edges which are not in any $K_{A_i,B_i}$  with $A_i \in \l_2$ is at least
\[
\left( \frac{p}{2}+o(1) \right)u^2-2p^2|\l_3|u-3p^3|\l_2 \setminus \l_3|u -8u^{3/2} \log^{1/2} n     - 8u^{3/2} \log^{3/2} n \log^{1/2} u.
\]
Observe that the expression above is a decreasing function if we view $|\l_3|$ as the  variable.
From the assumptions $|\l_3| <\left( \tfrac{1}{2}-\tfrac{1}{250} \right)u$  and $|\l_2| \leq \left(\tfrac{1}{2}+\tfrac{1}{1500}  \right)u$, the number of edges which are not contained in any $K_{A_i,B_i}$ with $A_i \in \l_2$ is at least
\[
\frac{p}{2}u^2-\left(\frac{1}{2}-\frac{1}{250}  \right)2p^2u^2-\frac{7}{500}p^3u^2+o(u^2) \geq \left(\frac{p^2}{125}-\frac{7p^3}{500} +o(1) \right) u^2
\]
when $n$ is large enough. Here we  note that $u^{3/2} \log^{3/2} n \log^{1/2} u=o(u^2)$ as we assume $u \geq c (\ln_b n)^{3+\epsilon}$.  Since $p \leq \tfrac{1}{2}$ and $p$ is a constant, we get that $\tfrac{p^2}{125}-\tfrac{7p^3}{500} $ is a positive constant.

Applying  Lemma \ref{l:lm33} with $\delta=\delta_2$,    the number of edges covered by   $ \l \setminus \l_2$ (i.e., $|A_i| \geq \delta_2 \ln_b n$) is at most
 \[
 \delta_2 |\l\setminus \l_2|u^{1-\delta_2} \ln_b u+8  u^{(3+\delta_2)/2} \log^{3/2} n \log^{1/2} u.
 \]
Since $u^{(3+\delta_2)/2} \log^{3/2} n \log^{1/2} u=o(u^2)$ by the  choice of $\delta_2$,
in order to cover the remaining edges,  we need at least  $C_1 u^{1+\delta_2/2}$ extra complete bipartite graphs for some positive constant $C_1$, i.e., $|\l \setminus \l_2 | \geq C_1 u^{1+\delta_2/2}$.  Since $C_1 u^{1+\delta_2/2}> \left( 1+\tfrac{1}{10000}\right)u$
for sufficiently large $n$,  this leads to a contradiction. Thus we have  $|\l_3| \geq (\tfrac{1}{2}-\tfrac{1}{250})$ and
 the claim is proved.

Now we proceed to prove the lemma using the fact that  $|\l_3| \geq \left(\tfrac{1}{2}-\tfrac{1}{250} \right)u$. We consider a auxiliary graph $U^{\ast}$  whose vertex set is $U$ and  edge set is $\l_3$.
 %It is possible that $U^{\ast}$ has multi-edges.
  We partition the vertex set of $U^{\ast}$ into three sets $U_1,U_2$, and $U_3$, where
$U_1=\{v \in V(U^{\ast}) \colon d_{U^{\ast}}(v)=0\}$, $U_2=\{v \in V(U^{\ast}) \colon d_{U^{\ast}}(v)=1\}$, and $U_3=\{v \in V(U^{\ast}) \colon d_{U^{\ast}}(v) \geq 2\}$. We will prove the following.

\noindent
{\bf Claim 2:}  The number of edges not contained in any $K_{A_i,B_i}$ for $A_i \in \l_3$ is at least $p\binom{|U_1|}{2}+ p^3 \binom{|U_2|/2}{2}+o(u^2)$.
\\
{\it Proof of Claim 2:}
 The first part of the sum follows from Lemma \ref{l:lm1}. For the second part of the sum, we let $U_2'$ be a maximum subset of  $U_2$ such that
 for each $v \in U_2'$, the neighbor of $v$ in $U^*$ is not in  $ U_2'$. We have $|U_2'| \geq |U_2|/2$.
 Then we  apply Lemma \ref{l:lm4} with $S=U_2'$ and   $T$ consisting of neighbors of $S'$ in $U^*$.
  To finish the proof of Claim 2, we
 use the fact  $u^{3/2} \ln n=o(u^2) $ as  $u \geq c(\ln_b n)^{3+\epsilon}$.

We will prove  Lemma  \ref{l:lm6} by contradiction. Suppose $|\l_3| \leq \left( \tfrac{1}{2}+\tfrac{1}{2000} \right)u$. This implies that
 the average degree of $U^{\ast}$ is at most $ 1+\tfrac{1}{1000}$. We consider the following cases.
\begin{description}
\item[Case 1:] $|U_3| \geq \tfrac{1}{5}u$.

  By considering the total sum of degrees of $U^{\ast}$,  we have
\[
2|U_3|+(u-|U_1|-|U_3|) \leq \left( 1+\frac{1}{1000} \right)u.
\]
 Thus, $|U_3| \geq \tfrac{1}{5}u$ implies $|U_1| \geq \tfrac{1}{6}u$. Claim 2 together with this  lower bound on $|U_1|$ implies that the number of edges not in any $K_{A_i,B_i}$ with $A_i \in \l_3$ is at least $\left(\tfrac{p}{72}+o(1) \right) u^2$.  By Claim 1 we have  $|\l_3| \geq \left( \tfrac{1}{2}-\tfrac{1}{250} \right)u$, so  the number of additional complete bipartite graphs $K_{A_i,B_i}$ with $A_i \in \l_2 \setminus \l_3$ is at most $\tfrac{7}{1500}u$ using the assumption $|\l_3| \leq \left( \tfrac{1}{2}+\tfrac{1}{2000} \right)u$. These complete bipartite graphs can cover at most $ (\tfrac{7p^3}{500}+o(1))u^2$ edges by Lemma \ref{l:lm32}. Thus we conclude  that the number of edges not covered by any of $K_{A_i,B_i}$ with $A_i \in \l_2$ is at least
\[
\left(\frac{p}{72}- \frac{7p^3}{500}+o(1) \right)u^2.
\]
Note that  $\tfrac{p}{72}- \tfrac{7p^3}{500}$ is a positive constant when $p$ is constant and $p \leq \tfrac{1}{2}$.
By applying   Lemma \ref{l:lm33} with $\delta=\delta_2$,   the bipartite graphs $K_{A_i,B_i}$ with  $ A_i \in  \l \setminus \l_2$ (i.e., $|A_i| \geq \delta_2 \ln_b n$) can cover at most
 \[
 \delta_2 |\l\setminus \l_2|u^{1-\delta_2} \ln_b u+8  u^{(3+\delta_2)/2} \ln^{3/2} n \log^{1/2} u
 \]
edges.  We note  $ u^{(3+\delta_2)/2} \ln^{3/2} n \log^{1/2} u=o(u^2)$ because of the  choice of $\delta_2$.  To cover the remaining edges,  we need at least  $C_1' u^{1+\delta_2/2}$ extra complete bipartite graphs $K_{A_i,B_i}$ with $A_i \in \l \setminus \l_2$ for some positive constant $C_1'$, i.e.,  $|\l \setminus \l_2 | \geq C_1' u^{1+\delta_2/2}$. Since $C_1' u^{1+\delta_2/2}> \left( 1+\tfrac{1}{10000}\right)u$ for $n$ large enough,
we get a  contradiction to the assumption  $|\l| \leq \left( 1+\tfrac{1}{10000}\right)u$.

\item[Case 2:] $|U_3| < \tfrac{1}{5}u$.

In this case we have $|U_1|+|U_2| \geq \tfrac{4}{5}u$. Note that the lower bound given by Claim 2 is minimized when $|U_2|=\tfrac{4}{5}u$, i.e., the number of edges not contained in any $K_{A_i,B_i}$ with $A_i \in \l_3$ is at least $\left( \tfrac{2}{25}p^3+o(1)  \right)u^2$. By the same argument as in Case 1 we can show the number of edges in $K_{A_i,B_i}$ with $A_i \in \l_2 \setminus \l_3$ is at most  $ (\tfrac{7p^3}{500}+o(1))u^2$.  Now there are at least
 \[
 \left( \frac{33}{500}p^3 +o(1) \right) u^2
 \]
edges which is not in any $K_{A_i,B_i}$ with $A_i \in \l_2$.  We note that  $\tfrac{33}{500}p^3$ is a positive constant  as we assume $p$ is a constant.  By using Lemma \ref{l:lm33} with $\delta=\delta_2$,  the  bipartite graphs $K_{A_i,B_i}$ with  $ A_i \in  \l \setminus \l_2$ (i.e., $|A_i| \geq \delta_2 \ln_b n$) can cover at most
 \[
 \delta_2 |\l\setminus \l_2|u^{1-\delta_2} \ln_b u+o(u^2)
 \]
 edges. As in Case 1, we  consider the number of extra bipartite graphs $K_{A_i,B_i}$ with $A_i \in \l \setminus \l_2$  needed to cover the remaining edges, leading to the same contradiction to the assumption on $\l$.
 \end{description}
 Therefore we have proved $|\l_3| > (\tfrac{1}{2}+\tfrac{1}{2000})$.
 \hfill $\square$

\noindent
{\bf Remark:}   When we apply Lemma \ref{l:lm33},    we require the error term is in a lower order of magnitude in comparison to the main term.  To make this satisfied, we have to
 assume $|U|=\Omega((\log n)^{3+\epsilon})$ .

\begin{lemma} \label{l:lm7}
Let $H$ be  a hypergraph  with the vertex set $U$ and the edge set $\l_1$.
There is some positive constant $C_2$ such that there are $C_2u$ vertices of $H$ with degree less than $\left( \tfrac{\delta_1}{2}-\tfrac{\delta_1}{3000}\right) \ln_b u$.
\end{lemma}

\noindent
{\bf Proof:} We consider several cases.

\begin{description}
\item [Case a:]   $|\l_2| > \left( \tfrac{1}{2}+\tfrac{1}{1500}\right) u$.
\\
The sum of  degrees in $H$ is less than
\begin{align*}
\delta_2 |\l_2| \ln_b u + \delta_1 |\l_1 \setminus \l_2| \ln_b u  &\leq \delta_2 \left( \frac{1}{2}+\frac{1}{1500}\right) u \ln_b u
+\delta_1 \left(\frac{1}{2}+\frac{1}{10000}-\frac{1}{1500}  \right) u \ln_b u\\
&\leq \left(\frac{\delta_1}{2} -\frac{\delta_1}{2000} \right) u \ln_b u.
\end{align*}
Here we used the assumption $|\l_1| \leq |\l|=m < \left( 1+\tfrac{1}{10000} \right)u$ and the choice of $\delta_2$.

\item [Case b:]
 $|\l_2| \leq \left( \tfrac{1}{2}+\tfrac{1}{1500}\right) u$.
 \\
By Lemma \ref{l:lm6},  $|\l_3| \geq \left( \tfrac{1}{2}+\tfrac{1}{2000} \right)u$. The sum of  degrees is at most
\begin{align*}
2|\l_3|+\delta_1 |\l_1\setminus \l_3| \ln_b u  &\leq 2 \left( \frac{1}{2}+\frac{1}{2000} \right)u+ \delta_1 \left( \frac{1}{2}+\frac{1}{10000}-\frac{1}{2000} \right) u \ln_b u \\
 &\leq \left( \frac{\delta_1}{2}-\frac{\delta_1}{2000} \right) u \ln_b u.
\end{align*}
\end{description}
We have proved that  the sum of  degrees of $H$ is less than $\left(\tfrac{\delta_1}{2} -\tfrac{\delta_1}{2000} \right) u \ln_b u$. Let $U'$ be the set of vertices with degree at least  $\left( \tfrac{\delta_1}{2}-\frac{\delta_1}{3000}\right) \ln_b u$. We consider
\[
|U'|\left( \frac{\delta_1}{2}-\frac{\delta_1}{3000}\right) \ln_b u \leq \left(\frac{\delta_1}{2} -\frac{\delta_1}{2000} \right) u \ln_b u,
\]
which yields $|U'| \leq (1-C_2)u$ for some positive constant $C_2$. Each vertex in $U \setminus U'$ has degree less than $\left( \tfrac{\delta_1}{2}-\tfrac{\delta_1}{3000}\right) \ln_b u$ and  $|U \setminus U'| \geq C_2u$. The lemma is proved. \hfill $\square$

We recall that $G[U]$ is the subgraph of $G$ induced by $U$. We have the following lemma.
\begin{lemma} \label{l:lm8}
The number of edges in $G[U]$ which are not contained in any $K_{A_i,B_i}$ with $A_i \in \l_1$ is at least $C_3 u^{2-\delta_1+\delta_1/2000}$ for some positive constant $C_3$.
\end{lemma}

\noindent
{\bf Proof:} We consider the hypergraph  $H$   with the vertex set $U$ and the edge set $\l_1$ as defined in Lemma \ref{l:lm7}.  Let $W$ be the set of vertices with degree less than $\left( \tfrac{\delta_1}{2}-\tfrac{\delta_1}{3000} \right)\ln_b u$ in $H$; we have $|W| \geq C_2 u$ for some positive constant $C_2$ by Lemma \ref{l:lm7}.

We will use Lemma \ref{l:lm5} to prove Lemma \ref{l:lm8}. In order to apply Lemma \ref{l:lm5}, we will first
 find a subset $W'$ of $W$ such that for any  $u,v \in W'$  there is no  $A_i \in \l_1$ containing $u$ and $v$.  Also  we will associate  each
 $w\in W'$  with a set $L(w) \subset U \setminus W'$  satisfying the property that $L(w) \cap L(w') =\emptyset$ for each $w \not =w' \in W'$.

To do so, we  consider an arbitrary linear ordering of vertices in $W$. Let $q= |W| / \ln_b^2 u$, $W_0=W$, $Z_0=\emptyset$ and $H_0=H$. For each $1 \leq i \leq q$, we recursively define a vertex $v_i$, a set $W_i$, a set $Z_i$ and a hypergraph $H_i$ as follows:  For given $W_{i-1}$ and $H_{i-1}$, we let $v_i$ be the first vertex in $W_{i-1}$ and define $F(v_i)=\{A \in E(H_{i-1}) \colon v_i \in A\}$.
 By the assumption on the size of sets in  $\l_1$ and  the degree upper bound  for vertices in  $W$, we have  $|\cup_{A \in F(v_i)} A| \leq \ln_b^2 u/2$. We define  $Z_i=\{A \in E(H_{i-1}) \colon |A  \setminus \left(\cup_{A' \in F(v_i)} A'\right)|=1\}$. Then
$|\cup_{A \in Z_i} A  \setminus \left(\cup_{A' \in F(v_i)} A'\right)| \leq \ln_b^2 u/2$ since each $A' \in F(v_i)$ can contribute at most $\delta_1 \ln_b u$ to the sum and $|F(v_i)| \leq \tfrac{\delta_1}{2} \ln_b u$ because of  the degree upper bound for vertices in   $W$.
We define $W_i=W_{i-1} \setminus \left(  \cup_{A \in Z_i \cup F(v_i)} A \right)$ and  $H_i$ to be the new hypergraph with the  vertex set $V(H_{i-1}) \setminus \left( \cup_{A \in Z_i \cup  F(v_i)} A \right)$. If $ A \in E(H_{i-1})$ then $A \setminus  \left( \cup_{A \in Z_i \cup F(v_i)} A \right) \in E(H_i)$.
  %We note if $w_j \in A$ for some $j \geq i$ and $A \in H_{i-1}$, then $A \setminus  \left( \cup_{A \in Z_ \cup F(v_i)} A \right)$ is still a set %containing $w_j$ (if it have size at least 2) and we have to keep it. This is the reason why we define the edge set of $H_i$ in this way.
Thus  $|V(H_i)| = |W_i| \geq |W_{i-1}|- \ln_b^2 u$ and  $W_{q-1} \not = \emptyset$. Therefore,  $v_i$ is well-defined for  $1 \leq i \leq q$. We write  $W'=\{v_1,v_2,\ldots,v_q\}$.

 For each $A \in F(v_i)$ and $A' \in  F(v_j)$ with $i<j$ we have $A \cap A'= \emptyset$ as we delete the set $\cup_{A \in Z_i \cup F(v_i)} A$ in step $i$. For each $v_i \in W'$ and each $A \in F(v_i)$,  we let $f(A)$ be an arbitrary vertex other than $v_i$ from $A$  and $F'(v_i)=\cup_{A \in F(v_i)} \{f(A)\}$. It follows from the preceding definitions that  $F'(v_i) \cap F'(v_j)=\emptyset$ for $1 \leq i \not =j \leq q$.  Now  an application of  Lemma \ref{l:lm5} with $U=V(H)$, $W=W'$,  $L(v_i)=F'(v_i)$ for each $v_i$  and $c= \tfrac{\delta_1}{2}-\tfrac{\delta_1}{3000}$ will prove the lemma.

We are left to verify the function $h(G,U,W,L)$ indeed gives a lower bound on the number of edges which is not covered by $\l_1$.
 From the construction , for each $v_i$ and each $A \in \l_1$ containing $v_i$,   either $A$ is  in  $F(v_i)$ or a subset of $A$ with size at least two is in $F(v_i)$. Hence,  $A \cap F'(v_i) \not =\emptyset$.  For an edge $\{v_i,v_j\}$, if
$\{v_i,z\}$ is a non-edge for each $z \in F'(v_j)$ and $\{v_j,z'\}$ is a non-edge for each $z' \in F'(v_i)$, then the edge $\{v_i,v_j\}$ is uncovered by the family of sets $\l_1$. Suppose $\{v_i,v_j\}$ is in $K_{A,B}$ for some $A \in \l_1$. We have either $v_i \in A$ or $v_j \in A$. In the former case we get   $ A \cap L(v_i) \not = \emptyset$ by the definition of $L(v_i)$. Let $z \in A \cap L(v_i)$.  Then $A$ and $B$ does not form a complete bipartite graph since $\{v_j,z\}$ is not an edge by the assumption. We get a contradiction and  we have a similar argument for the latter case.
\hfill $\square$

\noindent
{\bf Remark:} We mention here that when we defined the set   $W '  \subset W $,  we did not aim to find the largest one as $|W'|=|W|/\log_b^2 n$ is large enough for proving Theorem \ref{t:thm2}.

We are ready to prove Theorem \ref{t:thm2}.

\noindent
{\bf Proof of Theorem \ref{t:thm2}:}  Suppose that
\[
E(G[U])=\bigsqcup_{i=1}^m E(K_{A_i,B_i}).
\]
  If  $m > \left( 1+ \tfrac{1}{10000} \right)u$, then we are done.  Otherwise,  Lemma \ref{l:lm8} implies that there are at least $C_3u^{2-\delta_1+\delta_1/2000}$ edges uncovered after we delete the edges in $K_{A_i,B_i}$ for each $A_i \in \l_1$.
  We then apply    Lemma \ref{l:lm33} with $\delta=\delta_1$ which gives an upper bound for    the number of edges covered by $\l \setminus \l_1$ (i.e., $|A_i| \geq \delta_1 \ln_b u$)  :
   \[
   \delta_1|\l \setminus \l_1| u^{1-\delta_1} \ln_b u +8  u^{(3+\delta_1)/2} \ln^{3/2} n \ln^{1/2}  u.
   \]
Here we note that $u^{(3+\delta_1)/2} \ln^{3/2} n \ln^{1/2}  u=o(u^{2-\delta_1+\delta_1/2000})$ because of  the choice of $\delta_1$.
%\mnote{Here is the place where  the selection of $\delta_1$ comes from.}
   Therefore we need at least $C_4 u^{1+\delta_1/2500} $ additional complete bipartite graphs $K_{A_i,B_i}$  with $A_i \in \l \setminus \l_1$ to cover the remaining edges, where $C_4$ is some positive constant.  Since $C_4 u^{1+\delta_1/2500} > 1.0001 u$ when $n$ is sufficiently large and we get a contradiction. Theorem \ref{t:thm2} is proved. \hfill $\square$

\section{Proof of Theorem \ref{t:thm1}}\label{sec6}
Before proving Theorem \ref{t:thm1}, we first state the following  lemma. The  proof will be omitted as it is very simple.
\begin{lemma}
\label{decom}
Suppose that edges of $G$ can be decomposed  into $k_1$ complete bipartite graphs, of which $k_2$ complete bipartite graphs
are stars for some $k_2 \leq k_1$. Then $G$ has an edge decomposition
 %we relate $\\tau(G)$ to $\\tau'(G')$ for some induced subgraph $G'$ of $G$.
$E(G)=\sqcup_{i=1}^k E(K_{A_i,B_i})$  with $k \leq k_1$ such that for $i \leq k_2$,  $K_{A_i, B_i}$ are stars  and  for $j >k_2$,
we have $A_j, B_j \subseteq V(G)\setminus \cup_{i=1}^{k_2} A_i$.
\end{lemma}

We are ready to prove Theorem \ref{t:thm1}.

\noindent
{\bf Proof of Theorem \ref{t:thm1}:}   The upper bound follows from the well-known fact  (see Theorem \ref{alpha}) that asymptotically almost surely  a random graph  $G \in G(n,p)$ has an independent set $I$ with size $ (2+o(1))\ln_{1/(1-p)} n$ where   $p$ is constant.
We consider vertices $v_1,\ldots,v_m$ with $m=n-ги2+o(1))\ln_{1/(1-p)} n$, which are not contained in $I$ . For each $1 \leq i \leq m$  we define a star $K_{A_i,B_i}$ with $A_i=\{v_i\}$ and $B_i=\{v_j \colon j>i \textrm{ and } \{v_i,v_j\} \in E(G)\}$. We have
\[
E(G)=\bigsqcup_{i=1}^m E(K_{A_i,B_i}).
\]
Therefore we have $\tau(G) \leq n-(2+o(1))\ln_{1/(1-p)} n$.

For the lower bound, we may assume that    $G \in G(n,p)$  satisfies all statements in  the lemmas in the preceding sections. Suppose $G$ has  an edge decomposition:
\[
E(G)=\bigsqcup_{i=1}^k E(K_{A_i,B_i}),
\]
with $k=\tau(G)$ and  assume that  for some $l \leq k$, we have $A_i=\{v_i\}$ for $1 \leq i \leq l$.

Let $W=\{v_1, \ldots, v_l\}$.
If $W=\emptyset$ then  Theorem \ref{t:thm1} follows from Theorem \ref{t:thm2} directly.  We need only to consider the case $W \not =\emptyset$.
By Lemma \ref{decom},  we can assume $E(G')=\sqcup_{i=l+1}^k E(K_{A_i,B_i})$ where $G'$ is the subgraph induced by $T=V(G) \setminus W$. We get
\begin{equation}
\tau(G)=|W|+\tau'(G').
\end{equation}
We will prove   $l >  n-c (\ln_{1/p} n)^{3+\epsilon}$ for any positive constants $c$ and $\epsilon$. Suppose  $l \leq n-c (\ln_{1/p} n)^{3+\epsilon}$ for some $c$ and $\epsilon$.  Thus, $|T|\geq c (\ln_{1/p} n)^{3+\epsilon}$. By Theorem \ref{t:thm2} we have $\tau'(G[T]) \geq \left( 1+\tfrac{1}{10000} \right) |T|$. Therefore

  \begin{align*}
  \tau(G)=|W|+\tau'(G')  &\geq |W|+\left( 1+\tfrac{1}{10000} \right) |T|
                      \geq n,
 \end{align*}
 which is a contradiction.    Theorem \ref{t:thm1} is proved. \hfill $\square$

\section{Problems and remarks}\label{sec8}
%The conjecture of Erd\H{o}s,  that almost all graphs $G$ satisfies $\tau(G) = n-\alpha(G)$  as stated in (\ref{ep1}). Recently, Alon   (see Theorem 1.1 in \cite{alon}) disproved this conjecture by showing $\tau(G) \leq n-\alpha(G)-1$ for almost all graphs $G \in G(n,1/2)$. The following two variations of Erd\H{o}s'
The results on the bipartite decomposition  in this paper lead to many questions, several of which  we mention here.

 For  $G \in G(n,1/2)$,    the lower bound $\tau(G) \geq n-o((\log n)^{3+\epsilon})$ for any positive constant $\epsilon$ is given by Theorem \ref{t:thm0} in this paper.
For the upper bound,  the result from \cite{abh}  gives    asymptotically almost surely $  \tau(G) \leq n-(1+c) \alpha(G)$  for some positive constant $c$.  Similar upper bound is not known for  any constant $p$ with $p < 1/2 $.  We believe the following conjecture is true.

\vspace{.1in}

\noindent
{\bf Conjecture 1:}  For a random graph $G \in G(n,p)$,  where $p$ is a constant and $p < 1/2$,  asymptotically almost surely we have $\tau(G)=n-(2+o(1))\log_{1/(1-p)} n$.

%\vspace{.1in}
%
%\noindent
%{\bf Conjecture 1:}
%For a random graph $G \in G(n,p)$, where $p$ is a constant no greater than $1/2$,  asymptotically almost surely
%\[ \tau(G) = \min \{n-\alpha(G),n-\beta(G)+1\}, \]
%where $\beta(G)$ is the size of largest induced complete bipartite graph in $G$.
%
%\vspace{.1in}
%
%A somewhat weaker conjecture is the following:

%\noindent
%{\bf Conjecture 1:}
%For a random graph $G \in G(n,p)$, where $p$ is a constant less than $1/2$,   asymptotically almost surely
%\[ \tau(G) = n - (2+o(1)) \log_{1/(1-p)} n, \]

For sparser random graphs,
 Alon \cite{alon}   showed   that there exists some (small) constant $c$ such that for $\tfrac{2}{n} \leq p \leq c$, a random graph $G$ in $G(n,p)$ satisfies $\tau(G)=n-\Theta\left(\tfrac{\log np}{p}  \right)$. It will be of interest to further sharpen the lower bound.

\noindent
{\bf Conjecture 2:}
For a random graph $G \in G(n,p)$, with $ p = o(1)$,  asymptotically almost surely
\[ \tau(G) =n-(1+o(1))\tfrac{2}{p} \log np.\]

%
%\vspace{.1in}
%
%\noindent
%{\bf Conjecture 4:}
%For a random graph $G \in G(n,p)$, where $p$ is a constant no greater than $1/2$,   suppose that an edge decomposition
%$E(G) = \sqcup_{i=1}^k E(K_{A_i,B_i})$ achieves $\tau(G)=k$. Then asymptotically almost surely at least $n-o(n)$ of the bipartite subgraphs  $K_{A_i,B_i}$ are stars.

In this paper, we have given rather crude estimates for the constants involved.  In particular, for the strong bipartition number $\tau'(G)$, a consequence of Theorem \ref{t:thm2} states  that for $G \in G(n,p)$, where $p$ is a constant no greater than $1/2$, asymptotically almost surely $\tau'(G) \geq 1.0001 n$.
For the case of $p \leq c$ for some small $c$,
Alon  \cite{alon} showed  that asymptotically almost surely $\tau'(G) \geq 2 n$ for  $G \in G(n,p)$.
 A natural question is to improve the  lower bound for $\tau'(G)$.


\begin{thebibliography}{9}
%\bibitem{as} N.~Alon and J.~Spencer, {\it The Probabilistic Method,} 3rd ed., John Wiley \& Sons, Hoboken, NJ, 2008.
\bibitem{alon} N.~Alon, Bipartite decomposition of random graphs, {\it J.~Combin.~Theory Ser.~B}, {\bf 113} (2015), 220--235.
\bibitem{abh} N.~Alon, T.~Bohman, and H.~Huang, More on the bipartite decomposition of random graphs, arXiv: 1409.6165[math.CO].
\bibitem{azuma} K.~Azuma, Weighted sums of certain dependent random variables, {\it Tohuku Math.~Journal}, {\bf 19}(3) (1967), 357--367.


\bibitem{chernoff}
H.~Chernoff,
A note on an inequality involving the normal distribution,
{\it Ann.~Probab.}, {\bf 9} (1981),
533--535.
\bibitem{ep} P.~Erd\H{o}s, Some remarks  on the theory of graphs, {\it Bull. Amer. Math. Soc.}, {\bf 53} (1947), 292--294.
%\bibitem{frieze} A.~Frieze, On the independence number of random graphs, {\it Discrete Math.}, {\bf 81} (1990), 171--175.
\bibitem{gp1} R.~L.~Graham and H.~O.~Pollak, On the addressing problem for loop switching, {\it Bell
Syst.~Tech.~J.},  {\bf 50} (8) (1971), 2495--2519.
%\bibitem{gp2}  R.~L.~Graham and H.~O.~Pollak, On embedding graphs in squashed cubes, {\it in Graph
%theory and applications} pp. 99--110. Lecture Notes in Math., Vol.~303, Springer, Berlin, 1972.
\bibitem{gm} G.~R.~Grimmett and C.~J.~H.~McDirmid, On colouring random graphs, {\it Math. Proc. Cambridge Phil. Soc.}, {\bf 77} (1975), 313--324.
\bibitem{karp} R.~Karp, Reducibility among combinatorial problems, {\it Complexity of Computer Computations} (R.~E.~Miller, J.~W.~Thatcher eds.) Plenum, New York, (1972), 85--103.

\bibitem{krw}T.~Kratzke, B.~Reznick and D.~West, Eigensharp graphs: Decomposition into
complete bipartite subgraphs, {\it  Trans.~Amer.~Math.~Soc.~}, {\bf 308} (1988), 637--653.
%\bibitem{lp} L.~Lu and X.~Peng, Loose Laplacian spectra of random hypergraphs, {\it Random Structures and Algorithms}, {\bf 41}(4), (2012), %521--545.
\bibitem{lo} L.~Lov\'{a}sz, On coverings of graphs, {\it Theory of Graphs (Proc. Conf. Tihany)}, Academic Press, (1969), 231--236.
\bibitem{peck} G.~W.~Peck, A new proof of a theorem of Graham and Pollak, {\it Discrete Math.}, {\bf 49} (1984), 327--328.
\bibitem{tv} H.~Tverberg, On the decomposition of $K_n$ into complete bipartite graphs, {\it J.~Graph
Theory},  {\bf 6} (1982), 493--494.
\bibitem{v1} S.~Vishwanathan, A polynomial space proof of the Graham-Pollak theorem, {\it J.~Combin.~Theory Ser.~A}, {\bf 115} (2008), 674--676.
\bibitem{v2} S.~Vishwanathan, A counting proof of the Graham Pollak Theorem,  {\it Discrete Math.},  {\bf 313}(6) (2013), 765--766.
\bibitem{yan}  W.~Yan and Y.~Ye, A simple proof of Graham and Pollak's theorem,  {\it J.~Combin.~Theory Ser.~A,} {\bf  113}, (2006), 892--893.

\end{thebibliography}
\end{document}